\documentstyle[11pt]{article}
\newtheorem{theorem}{Theorem}[section]
\newtheorem{definition}[theorem]{Definition}
\newtheorem{proposition}[theorem]{Proposition}
\newtheorem{corollary}[theorem]{Corollary}
\newtheorem{lemma}[theorem]{Lemma}

\def \proof {\noindent {\bf Proof.}\ \ }
\def \endproof {{\mbox{}\nolinebreak\hfill\rule{2mm}{2mm}\par\medbreak} }

\def \remarks {\noindent {\bf Remarks.}\ \ }

\newcommand {\appr}[1] {\stackrel{#1}{\approx}}
\newcommand {\text}[1] {\mbox{\ \ \ #1 \ \ }}

\def \Pelc {A.~Pelczy\'{n}ski}

\def \e {\varepsilon}
\def \d {\delta}
\def \l {\lambda}
\def \O {\Omega}
\def \G {\Gamma}
\def \f {\varphi}
\def \F {\Phi}
\def \G {\Gamma}
\def \N {{\bf N}}
\def \R {{\bf R}}

\def \# {\vspace{0.5cm} \noindent}

\def \ra {\rightarrow}
\def \span {{\rm span}}
\def \supp {{\rm supp}}
\def \free {{\rm free}}
\def \dist {{\rm dist}}

\def \nm {{n \le m}}
\def \npm {{n \le p(m)}}
\def \nqm {{n \le q(m)}}
\def \nrm {{n \le r(m)}}
\def \ng {{n \ge 1}}
\def \mg {{m \ge 1}}
\def \jg {{j \ge 1}}
\def \pm {{p(m)}}
\def \qm {{q(m)}}
\def \nze {{n_0}}

\def \xnxn {(x_n, x_n^*)_{n \ge 1}}

\def \znzn {(z_n, z_n^*)_{n \ge 1}}
\def \xns {x_n^*}
\def \yns {y_n^*}
\def \zns {z_n^*}

\def \xh    {\widehat{x}}

\def \ehn {\widehat{e}_n}
\def \ehi {\widehat{e}_i}
\def \pin {{\pi(n)}}
\def \ymn {y_{m,n}}
\def \ymns {y_{m,n}^*}
\def \ymks {y_{m,k}^*}
\def \loo {l_2^{|\O(r(m))|}}

\begin{document}
\title {On constructions of strong and uniformly minimal M-bases
        in Banach spaces}
\author {R. Vershynin\thanks{This research was finished when the author
                             was visiting the Politecnic institute of Milan}}
\date{\scriptsize 1991 Mathematics Subject Classification: 46B15.}
\maketitle

\begin{abstract}
  We find a natural class of transformations ("flattened perturbations")
  of a norming M-basis in a Banach space $X$, 
  which give a strong norming M-basis in $X$.
  This simplifies and generalizes the positive answer
  to the "strong M-basis problem" solved by P.~Terenzi.
  We also show that in general one cannot achieve
  uniformly minimality applying standard transformations
  to a given norming M-basis, despite of the existence in $X$
  a uniformly minimal strong M-bases.
\end{abstract}

\section{Introduction}

Does every separable Banach space $X$ have a strong M-basis? 
This problem remained open for a long time (see \cite{S}, Problem 8.1)
and was solved in positive by P.~Terenzi (\cite{T 90}, \cite{T 94}).
He proved that every complete norming biorthogonal system 
has a block perturbation which is a strong complete norming 
biorthogonal system
(a biorthogonal system $\znzn$ is a {\em block perturbation} of 
 a biorthogonal system $\xnxn$ if for every $\mg$
\begin{equation}                                              \label{perturb}
  [z_n]_{n \in I(m)}   = [x_n]_{n \in I(m)}  \text{and}
  [z_n^*]_{n \in I(m)} = [x_n^*]_{n \in I(m)},
\end{equation}
where $I(m)$, $m=1,2,\ldots$, are some successive intervals 
of positive integers).

In Section \ref{strongbp} the way of constructing strong block perturbations
is essentially simplified and slightly generalized.
For a given complete norming biorthogonal system we find a certain class 
of block perturbations
(so-called "flattened perturbations") 
which are {\em strong} norming complete biorthogonal systems.
Therefore, we demonstrate a new construction of strong M-bases in every
separable Banach space. The presentation of this part is self-contained.

The second part of the paper is concerned with 
questions of uniform minimality.
A biorthogonal system $\znzn$ is called a {\em pile perturbation}
of a biorthogonal system $\xnxn$ if (\ref{perturb}) holds 
for $m=1,2,\ldots$, where $I(m)$ are some intervals of integers
with the left bounds $=1$ 
and the right bounds $\ra \infty$ as $m \ra \infty$.
The notions of block  and pile perturbations can be considered
for minimal sequences as well:
$(z_n)_\ng$ is called a block  (resp. pile) perturbation of 
$(x_n)_\ng$ if we have 
\begin{equation}                                             \label{perturb2}
  [z_n]_{n \in I(m)}   = [x_n]_{n \in I(m)}                   
\end{equation}
instead of (\ref{perturb}) in the correspondent definitions.
Clearly, each block perturbation is a pile perturbation.
It was a long standing open problem in Banach space theory
whether every separable Banach space $X$ has 
a uniformly minimal M-basis. R.~Ovsepian and \Pelc\ 
solved it in positive (\cite{OP}, see also \cite{P}).
Later P. Terenzi constructed in every $X$ 
a strong uniformly minimal M-basis (\cite{T 90}, 
see \cite{T 98} for further improvements).
In view of the results of Section \ref{strong},
these observations suggest the following question:
can one achieve uniform minimality in constructions of
strong block (or, at least, pile) perturbations?

The answer is in general negative.
In Section \ref{unimiR} we construct a 
complete norming biorthogonal system $\xnxn$
in $l_2$ without uniformly minimal pile perturbations.
It follows that the norming M-basis $(x_n)_\ng$ 
has no uniformly mininal block perturbations.
The system $\xnxn$ has even more pathological structure,
to be established in Theorem \ref{unb}.


I am grateful to V.~Kadets for the guidance, 
and for P. Terenzi for his hospitality 
during my visit to Milan.

\section{Standard definitions}

Usual preliminaries can be found in \cite{LT} and \cite{S}.
Nevertheless we recall some definitions.
A system $\xnxn \subset X \times X^*$ is called {\em biorthogonal}
if $\xns(x_m) = \d_{n,m}$ for every $n,m$ (Kronecker's delta).
Suppose a complete sequence $(x_n)_\ng$ is {\em minimal},
i.e. $x_n \not\in [x_m]_{m \ne n}$ for any $n$.
Then there exists a unique sequence of 
{\em biorthogonal functionals} $(\xns)_\ng \subset X^*$,
i.e. such that $\xnxn$ is a biorthogonal system.
A sequence $(\xns)_\ng \subset X^*$ is called {\em total} 
if for every $x^* \in X^*$ one can find an $n$ so that $\xns(x) \ne 0$.
Further, $(\xns)_\ng$ is called {\em norming}  
if there exists a constant $c>0$ such that 
for every $x^* \in X^*$ one can find an $x \in [x_n]_\ng$ 
with $|\xns(x)| \ge c\|x^*\|\|x\|$.
Trivially every norming sequence is total.
We will call a biorthogonal system itself $\xnxn$ complete 
if the sequence $(x_n)_\ng$ is complete, 
and total (resp. norming) 
if the sequence $(\xns)_\ng$ is total (resp. norming). 

A complete minimal system $(x_n)_\ng$ is called an {\em M-basis}
(resp. {\em norming M-basis}) if its sequence 
of biorthogonal functionals is total (resp. norming).
A complete total biorthogonal system $\xnxn$ 
(or simply an M-basis $(x_n)_\ng$) is called {\em strong}
if $x \in [\xns(x)x_n]_\ng$ for every $x \in X$.
There is an intrinsic characterization of strongness,
due to A.~Plans and A.~Reyes \cite{PR}:
an M-basis $(x_n)_\ng$ is strong iff 
$[x_n]_{n \in A} \cap [x_n]_{n \in B} = [x_n]_{n \in A \cap B}$
for every subsets of indices $A$ and $B$.

We say that a system $\xnxn \subset X \times X^*$ 
(not necessarily biorthogonal) is {\em $C$-bounded}
if $\|x_n\|\|\xns\| \le C$ for every $n$. 
Clearly, a complete biorthogonal system $\xnxn$ is $C$-bounded for some $C>0$
iff the sequence $(x_n)_\ng$ is {\em uniformly minimal},
i.e. $\inf_n \dist(x_n/\|x_n\|, [x_m]_{m \ne n}) > 0$ for every $n$.
In this case we call the system $\xnxn$ itself uniformly minimal.

\section{Strong block perturbations}                         \label{strongbp}

A partition of $\N$ into finite sets $(A(j))_\jg$
is called a {\em block partition} if for some successive 
intervals of integers $(I(m))_\mg$ the sets 
$\cup_{j \in I(m)} A(j)$
are successive intervals of integers, $m=1,2, \ldots$

We shall use the notion of block perturbations also for finite systems:
$(z_n, \zns)_{n \le m}$ is a block perturbation of $(x_n, \xns)_{n \le m}$
if $[z_n]_{n \le m}  = [x_n]_{n \le m}$  and 
   $[\zns]_{n \le m} = [\xns]_{n \le m}$.

Let $\xnxn$ be a biorthogonal system. 
Fix some partition of $\N$ into finite sets $(A(j))_\jg$
and a sequence of numbers $n(j) \in A(j)$.

\begin{definition}                                               \label{flat}
  A biorthogonal system $\znzn$ is called 
  a {\em flattened perturbation} of a biorthogonal system $\xnxn$ 
  with respect to $(n(j),A(j))$
  if for every $\jg$

  (i) $(z_n, \zns)_{n \in A(j)}$ is a block perturbatoin of  
      $(x_n, \xns)_{n \in A(j)}$;

  (ii) $\| \zns - x_{n(j)}^* \|  \le  \e_j / \|x_{n(j)}\|$ for $n \in A(j)$,
  
  \noindent where $\e_j$ are some positive scalars with $\sum \e_j < \infty$.
\end{definition}

Trivially, if $(A(j))$ is a block partition, 
then every flattened perturbation with respect to $(n(j),A(j))$
is a block perturbation.
Note that flattened perturbations are easy to construct:
one can apply an invertible linear operator acting in $[\xns]_{n \in A(j)}$
which sends each $\xns$ to some vector close to $x_{n(j)}^*$.

Now we state the main result in this section.

\begin{theorem}                                                \label{strong}
  Let $\xnxn$ be a complete norming biorthogonal system.
  Then there is a block partition $(A(j))$ 
  and numbers $n(j) \in A(j)$ so that 
  each flattened perturbation of $\xnxn$ with respect to $(n(j),A(j))$
  is a strong complete biorthogonal system.
\end{theorem}

We wil use the following two known results, due to P.~Terenzi.
Since their proofs are scattered among different papers, 
and for the sake of completeness, we prove these results below.

\begin{lemma}                                                    \label{repr}
  Let $\xnxn$ be a complete biorthogonal system in a Banach space $X$.
  Then there is a sequence of positive integers $r(1) < r(2) < \ldots$
  (which we call {\em representing indices})
  so that for every $x \in X$
  $$
  x = \lim_m \Big( \sum_{n=1}^{r(m)} \xns(x)x_n  + v_m \Big)
  $$
  for some vectors $v_m \in [x_n]_{n=r(m)+1}^{r(m+1)}$ depending on $x$.
\end{lemma}

\begin{lemma}                                                  \label{series}
  If $\xnxn$ is a complete norming biorthogonal system in $X$, 
  then representing indices can be chosen with the following property.
  Suppose for some $x \in X$ there is a sequence of
  positive integers $m_1 < m_2 < \ldots$ so that
  $$
  \text{the series}
  \sum_{k=1}^\infty \sum_{n=r(m_k)+1}^{r(m_k+1)} \xns(x)x_n 
  \text{converges.}
  $$
  Then setting $r(m_0)=0$ we have
  \begin{equation}                                               \label{xsum}
    x = \sum_{k=0}^\infty \sum_{n=r(m_k)+1}^{r(m_{k+1})} \xns(x)x_n.
  \end{equation}
\end{lemma}

In the sequel, the relation $x \appr{\e} y$ between two vectors 
$x$ and $y$ means that $\|x-y\| \le \e$. 
We also assume for convenience that $\sum_{n \in \emptyset} y_n = 0$
for every vectors $y_n$.

\# {\bf Proof of Lemma \ref{repr}.}
We proceed by induction. 
Set $r(1)=1$ and assume that $(r(n))_{n \le m}$ is constructed 
for some $\mg$. 
Then, by a simple compactness argument (hint: a finite net)
there exists a number $p(m+1) > r(m)$ large enough 
so that for every $z \in [x_n]_{n \le r(m)}$ with $\|z\|=1$ 
\begin{equation}                                                 \label{dist}
  \dist(z,[x_n]_{n=r(m)+1}^{p(m+1)})  \appr{\d}
  \dist(z,[x_n]_{n=r(m)+1}^\infty),  
\end{equation}
where $\d$ can be taken sufficiently small: 
$\d = ( m \sum_{n \le r(m)} \|\xns\|\|x_n\| )^{-1}$. 
Now set $r(m+1):=p(m+1)$.

Let us check that $(r(n))$ is, indeed, a sequence of representing indices. 
Let $x\in B(X)$ and $\e \in (0,1)$. 
If a number $m$ is sufficiently large, 
then there is an $\xh \in [x_n]_{n \le r(m)}$ with $\xh \appr{\e} x$. 
Then (\ref{dist}) holds for 
$z := \xh - \sum_{n \le r(m)} \xns(x)x_n$   with   $\d =\e$.
But $z \appr{\e} x - \sum_{n \le r(m)} \xns(x)x_n  =: x'$; 
hence 
$\dist( x', [x_n]_{n=r(m)+1}^{p(m+1)} )  \appr{3\e}
 \dist( x', [x_n]_{n=r(m)+1}^\infty )  = 0. $
Then $x'$ is within a distance $3\e$ from 
$[x_n]_{n=r(m)+1}^{p(m+1)} = [x_n]_{n=r(m)+1}^{r(m+1)}$.
This completes the proof.            \endproof

\# {\bf Proof of Lemma \ref{series}.}
Assume $(r(n))_{n \le m}$ is constructed for some $m\ge 0$. 
As in the proof of Lemma \ref{repr}, 
we find a number $p(m+1) > r(m)$ so that (\ref{dist}) holds 
for every $z \in [x_n]_{n \le r(m)}$ with $\|z\|=1$.
Recall that the biorthogonal system $\xnxn$ is $(2c)$-norming for some $c>0$. 
That is, given a $v \in X$, we have 
$x^*(v) \ge 2c\|v\|$ for some $x^* \in [\xns]_\ng$ with $\|\xns\|=1$. 
The simple compactness argument provides a number 
$r(m+1) > p(m+1)$ such that we have the following property:
\begin{description}          
  \item[(P)]
    If $v\in [x_n]_{n \le p(m+1)}$,
    then $x^*(v) \ge c\|v\|$ for some $x^* \in [\xns]_{n \le r(m+1)}$
    with $\|x^*\| =1$.  
\end{description}
Then $(r(n))$ constructed in this way is, indeed, 
a sequence of representing indices. 
Let us verify that the conclusion of Lemma \ref{series} holds.
Substracting the convergent series, we can assume that for every $k \ge 1$
\begin{equation}                                               \label{iszero}
  \xns(x)=0   \text{whenever} r(m_k)+1 \le n \le r(m_k+1).
\end{equation}
We can also assume that $x \in B(X)$.
Let $\e>0$; by the proof of Lemma \ref{repr} we have 
for any sufficiently large integer $k$
$$
x \appr{\e}  \sum_{n \le r(m_k)} \xns(x)x_n + v_{m_k}   \text{with}
v_{m_k} \in [x_n]_{n=r(m_k)+1}^{p(m_k+1)}.  
$$
Therefore, to finish the proof it is enough to show that 
$\lim_k v_{m_k} = 0$.

Using (\ref{iszero}), we get  
$v_{m_k}  \appr{\e}  x - \sum_{n \le r(m_k+1)} \xns(x)x_n$. 
Therefore $x^*(v_{m_k}) \appr{\e} 0$ 
for any $x^* \in [\xns]_{n \le r(m_k+1)}$ with $\|x^*\|=1$. 
On the other hand, by the property (P), we have
$x^*(v_{m_k}) \ge c\|v_{m_k}\|$  
for some $x^* \in [\xns]_{n \le r(m_k+1)}$ with $\|x^*\|=1$. 
These estimates yield $c\|v_{m_k}\| \le \e$. 
This completes the proof.                     \endproof

\# {\bf Proof of Theorem \ref{strong}.}
Let $(r(m))$ be representing indices of $\xnxn$.
We construct the block partition $(A(j))$ and numbers $n(j) \in A(j)$
by induction. At each successive step, we find a successive interval
of integers ending at some representing index $r(m)$, 
and this interval will be the new block of sets $A(j)$.
Suppose $(n(j),A(j))_{j \le j_0}$ is constructed 
and let $r(m)$ be the last element of the interval 
$\cup_{j \le j_0} A(j)$. We call such $r(m)$ a {\em block bound}.

For $r(m)+1 \le j \le r(m+1)$, let $d_j = m+j-r(m)$.
Let $E(j)$ be the set consisting of $\{j\}$ plus the successive interval
between representing indices:
$$
E(j) = \{j\} \cup \{ r(d_j)+1, \ldots, r(d_j+1) \}.
$$
We see that the sets $E(j)$ are disjoint and their union is an interval
beginning at $r(m)+1$.
These $E(j)$ will form the new block of sets $A(j)$. 
More precisely, we define
\begin{equation}
  A(j_0+j-r(m)) = E(j)  \text{and}   n(j_0+j-r(m)) = j           \label{ajnj}
\end{equation}
for $r(m)+1 \le j \le r(m+1)$. This completes the construction.

Observe that the set $\{n(1), n(2), \ldots\}$ is exactly the union of 
the intervals $\{r(m)+1, \ldots, r(m+1)\}$, where $r(m)$ are the
block bounds.
Let $j(n)$ be the (one-to-one) function from this set to $\N$ 
which maps $n(j)$ to $j$.

Now we verify the conclusion of Theorem \ref{strong}.
Let $\znzn$ be any flattened perturbation of $\xnxn$
with respect to $(n(j),A(j))$,
and let $\sum \e_n$ be a convergent series of positive numbers.
Pick any $x \in S(X)$. We are to show that $x \in [\zns(x)z_n]_\ng$.
There are two possibilities:
\begin{description}
 \item[(A)] There exists a sequence of block bounds 
   $r(m_1) < r(m_2) < \ldots$ such that for every $k \ge 1$
   $$
   \|\xns(x)x_n\| \le \e_{j(n)}  
   $$
   for $r(m_k)+1 \le n \le r(m_k+1)$;
 \item[(B)] For every sufficiently large block bound $r(m)$
   there is a $\nze=\nze(m)$ with $r(m)+1 \le \nze \le r(m+1)$ such that
   \begin{equation}                                                 \label{B}
     \|x_\nze^*(x)x_\nze\| > \e_{j(\nze)}  \text{and}
     \|\xns(x)x_n\| \le \e_{j(n)}  
   \end{equation}
   for $n_0 \le n \le r(m+1)$.
\end{description}
If (A) is the case, then the convergence of the series 
$ \sum_{k=1}^\infty \sum_{n=r(m_k)+1}^{r(m_k+1)} \e_{j(n)}  \le
\sum_{j=1}^\infty \e_j $
makes possible to apply Lemma \ref{series}.
We derive from (\ref{xsum}) that $x \in [\zns(x)z_n]_\ng$,
since $(z_n,\zns)_{n=r(m_k)+1}^{r(m_{k+1})}$ is a block perturbation of
$(x_n,\xns)_{n=r(m_k)+1}^{r(m_{k+1})}$.
This completes the proof in this case.

If (B) is the case, then another argument works.
Consider the set $\O = \{n \in \N : \zns(x)=0\}$.  
It is enough to show that $x \in [z_n]_{n \in \N \setminus \O}$.
Fix an $\e>0$ and a sufficiently large block bound $r(m)$;
let $n_0$ be an index guaranteed by (B).

\noindent CLAIM: \ $E(\nze)) \subset \N \setminus \O$.

Indeed, it follows from (\ref{ajnj}) that $E(n_0)=A(j(\nze))$.
If the Claim were not true, then $\zns(x)=0$ for some $n \in A(j(\nze))$.
Then by the definition of a flattened perturbation we would have 
$|x_\nze^*(x)|  \le  \e_{j(\nze)} / \|x_\nze\|$,  
which contradicts to (\ref{B}).

By Lemma \ref{repr}, there is a vector
$v \in [x_n]_{n \in E(\nze)}$  such that setting 
$\G = \{1, \ldots, r(m)\} \cup \big( \cup_{j=r(m)+1}^{\nze-1} E(j) \big)$
we have
\begin{eqnarray*}
  x  &\appr{\e}&  \sum_{n=1}^{r(d_\nze)} \xns(x)x_n + v                 \\
     &    =    &  \sum_{n \in \G} \xns(x)x_n 
                   + \big( x_\nze^*(x)x_\nze + v \big)
                   + \sum_{n=\nze+1}^{r(m+1)} \xns(x)x_n.
\end{eqnarray*}
The first summand belongs to $[z_n]_{n \in \N \setminus \O}$.
Indeed, $(x_n,\xns)_{n \in \G}$ is a block perturbation 
of $(z_n,\zns)_{n \in \G}$; thus 
$\sum_{n \in \G} \xns(x)x_n = \sum_{n \in \G} \zns(x)z_n
 \in [z_n]_{n \in \N \setminus \O}$.
The second summand belongs to 
$[x_n]_{n \in E(\nze)} = [z_n]_{n \in E(\nze)} 
 \subset [z_n]_{n \in \N \setminus \O}$  by the Claim.
The third summand has the norm less than $\e$ 
if $m$ and, therefore, $\nze=\nze(m)$, were chosen sufficiently large: 
this follows from (\ref{B}).

Thus we have shown that $\dist(x,[x_n]_{n \in \N \setminus \O}) < 2\e$.
This completes the proof.                    \endproof

\section {Uniformly minimal pile perturbations}     \label{unimiR}

We shall find a biorthogonal system $\xnxn$ in $l_2$
which has no uniformly minimal pile perturbations. 
This will follow from a more general result.

\begin{definition}
  We say that a biorthogonal system $\znzn$ {\em is spanned} by
  a biorthogonal system $\xnxn$ if
  \begin{eqnarray}                                            \label{spanned} 
    (z_n)_\ng \subset \span(x_n)_\ng  \text{and}     
    (\zns)_\ng \subset \span(\xns)_\ng.
  \end{eqnarray}
\end{definition}

\noindent If (\ref{spanned}) holds, 
then to every positive integer $\mg$ we can assign 
the minimal number $\qm$ such that 
$$
[z_n]_\nm \subset [x_n]_\nqm   \text{and}     
[\zns]_\nm \subset [\xns]_\nqm.
$$
We call $(q(m))_\mg$ {\em the spanning indices}.
Obviously, $\qm \ge m$ for every $\mg$.
Clearly, $\znzn$ is a pile perturbation of $\xnxn$ iff the equality
$\qm = m$ holds for infintely many positive integers $m$.

The main result in this section states that there are 
complete norming biorthogonal systems in $l_2$
such that uniformly minimal systems spanned by them must have
very large spanning indices.

\begin{theorem}                                                   \label{unb}
  Given a sequence of positive numbers $(\l_m)_\mg$, 
  there is a complete biorthogonal system $\xnxn$ in $l_2$ 
  with the following property.
  If $(\qm)$ are the spanning indices 
  of a uniformly minimal system spanned by $\xnxn$, then 
  $$
  \lim_m  \qm / \l_m  =  \infty. 
  $$
\end{theorem}

\begin{corollary}                                             \label{without}
  There exists a complete norming biorthogonal system 
  $\xnxn$ in $l_2$ 
  without uniformly minimal pile and block perturbations.
  Moreover, the norming M-basis $(x_n)_\ng$ has no
  uniformly minimal block perturbations.
\end{corollary}

\proof  The first statement follows from Theorem \ref{unb}
if we set $\l_m = m$, $m=1,2,\ldots$

Let $(z_n)_\ng$ be arbitrary block perturbation of $(x_n)_\ng$.
Let $I(m)$ be successive intervals of integers
so that (\ref{perturb2}) holds for every $m$. Then
$$
[\xns]_{n \in I(m)} 
 = \big( [x_n]_{n \in I(m'), m' \ne m} \big)^\perp
 = \big( [z_n]_{n \in I(m'), m' \ne m} \big)^\perp
 = [\zns]_{n \in I(m)}.
$$
Hence $\xnxn$ is a block perturbation of $\znzn$. 
Then it follows from the first part that $(x_n)_\ng$ 
is not uniformly minimal.  \endproof

\remarks   1. Of course, $(x_n)_\ng$ has a uniformly minimal 
pile perturbation (apply the standard biorthogonalization 
procedure in $l_2$).

2. It will follow from the proof that these results hold not only in $l_2$, 
but also in every reflexive Banach space with unconditional basis.

\#
We proceed now to the proof of Theorem \ref{unb}.

Let $(e_n)_\ng$ denote the canonical basis in $l_2$.
Let $\pi$ be some permutation on $\N$. We shall specify $\pi$ later,
it will depend only on the sequence $(\l_n)_\ng$.
Let $(\e_n)_\ng$ be a sequence of positive numbers such that
\begin{eqnarray}
  \sum_{i>n} \e_i^2  \le  1/8.                                    \label{eps}
\end{eqnarray}

The following proposition can be derived easily from the standard 
construction of an M-basis in a separable Banach space 
(see \cite{LT}, Proposition 1.f.3).

\begin{proposition}                                                \label{xn}
  There is a biorthogonal system $\xnxn$ in $l_2$ 
  and a sequence $(\ehn)_\ng$ in $l_2$ such that:

  (i)  $[x_n]_\nm = [\ehn]_\nm$  and  $[\xns]_\nm = [e_\pin]_\nm$  
       for every $\mg$;

  (ii) $\ehn \appr{\e_n} e_n$  for every $\ng$.
\end{proposition}

\noindent Apply Proposition \ref{xn} 
and define a linear operator $T$ in $l_2$ by
$$
T\ehn=e_n, \;\;\; \ng.
$$
It is not hard to check that 
(\ref{eps}) yields that $T$ is well defined and is an isomorphism:
\begin{equation}                                                   \label{tt}
  \|T\| \le 2,  \;\;\;  \|T^{-1}\| \le 2.
\end{equation}
In particular, $(\ehn)_\ng$ is a basis in $l_2$. 
Then $(x_n)_\ng$ and $(\xns)_\ng$ are complete sequences. 
Hence $\xnxn$ is a complete norming biorthogonal system.

It will {\em not} be enough to know that $T$ is just an isomorphism.
The following lemma shows that $T$ is asymptotically close to the identity.

\begin{lemma}                                                  \label{tozero}
  Let $(z_n)_\ng$ be a normalized M-basis in $l_2$. Then
  $$
  \lim_n (Tz_n-z_n) = 0.
  $$
\end{lemma}

\proof     It suffices to show that $\lim_n \|z_n - T^{-1}z_n\| = 0.$
Let $\e>0$. Let $n_0=n_0(\e)$ be a positive integer which we specify below. 
For every $\ng$, write the expansions
$z_n = \sum_i a_{n,i}e_i$,
where $(a_{n,i})$ are some scalars.
By the triangle inequality, H\"{o}lder's inequality,
and our choise of vectors $\ehi$, 
we have for any $k \ge 1$
\begin{eqnarray}                                                  \label{dif}
   \|z_n - T^{-1}z_n\| 
   & = &  \Big\| \sum_{i\ge 1} a_{n,i}(e_i-\ehi) \Big\|   \nonumber \\
   &\le&  \sum_{i \le k} |a_{n,i}| + 
          \Big( \sum_{i>k} |a_{n,i}|^2 \Big)^{1/2} 
          \Big( \sum_{i>k} \|e_i-\ehi\|^2 \Big)^{1/2}   \nonumber \\ 
   &\le&  \sum_{i \le k} |a_{n,i}| + 
          \Big( \sum_{i>k} \e_i^2 \Big)^{1/2}      
\end{eqnarray}

\noindent CLAIM: \  $\lim_n a_{n,i} = 0$ for every $i \ge 1$.

Indeed, since $(z_n)_\ng$ is a normalized M-basis in a reflexive space,  
the sequence $(z_n)_\ng$ tends weakly to zero. 
Then for every $i \ge 1$
$$
0 = \lim_n \langle  e_i, z_n  \rangle = \lim_n a_{n,i}.
$$
This proves the Claim.

Now we describe how to pick $n_0$.
First choose $k=k(\e)$ so that 
the second summand in (\ref{dif}) is less then $\e/2$.
By Claim, we can pick $n_0=n_0(k,\e)$ so that         
the first summand in (\ref{dif}) is less then $\e/2$  whenever  $n>n_0$.
Thus $\|z_n - T^{-1}z_n\| < \e/2+\e/2 = \e$
for every $n>n_0$.         \endproof

We proceed to verification of the conclusion of Theorem \ref{unb}.
Let $\znzn$ be a system spanned by $\xnxn$, and such that
\begin{equation}                                              \label{bounded}
  \|z_n\| = 1, \;\;\;  \|\zns\| \le M  \text{for} \ng.
\end{equation}
Assume that the conclusion of Theorem \ref{unb} is false.
Then there are a positive (integer) constant $c$
and increasing sequences of positive integers $(\pm)_\mg$ and $(\qm)_\mg$ 
such that for every $\mg$ we have:
\begin{eqnarray}
 r(m):=q(\pm) \le c\l_\pm                                        \label{qmpm}
\end{eqnarray}
and  
\begin{equation}
  [z_n]_\npm    \subset   [x_n]_\nrm = [\ehn]_\nrm,             \label{supp1}  
\end{equation}
\begin{equation}  
  [\zns]_\npm   \subset   [\xns]_\nrm = [e_\pin]_\nrm.          \label{supp2}
\end{equation}
For a positive integer $k$, set 
$$
\O(k) = \{1,2, \ldots, k\}  \cap  \{\pi(1),\pi(2), \ldots, \pi(k)\}.
$$
Note that if $k$ is large enough, then $\O(k)$ is not empty.

For any integer $m$ large enough, let $P_m$ denote the orthogonal
projection in $l_2$ onto $[e_n]_{n \in \O(r(m))}$. 
Define a system $(\ymn, \ymns)_\npm$ by
$$
\ymn = P_m T z_n,  \;\;\;\;\;  \ymns = P_m \zns.
$$
By (\ref{tt}) and (\ref{bounded}), this system is $(2M)$-bounded. 
One can not assert that it is biorthogonal, but this is not far from truth.

\begin{definition}
  Let $X$ be a Banach space and $\e>0$. 
  A system $(y_n, \yns) \subset X \times X^*$ (finite or infinite) 
  is called {\em $\e$-roughly biorthogonal} in $X$ if, 
  for every indices $k$ and $n$,
  $y_k^*(y_n) \appr{\e} \d_{k,n}$.
\end{definition}

\begin{lemma}
  There are positive integers $n_0$ and $m_0$ such that, 
  for every $m>m_0$, the system $(\ymn,\ymns)_{n=n_0+1}^\pm$
  is a $(1/4)$-roughly biorthogonal system in $l_2$.
\end{lemma}

\proof     Lemma \ref{tozero} provides a number $n_0$ such that    
\begin{equation}                                               \label{tzn-zn}
  \|Tz_n-z_n\| < \frac{1}{4M}  \text{for}  n>n_0.   
\end{equation}
Pick $m_0$ so that $\O(r(m_0))$ is not empty and $p(m_0)>n_0$.
Fix a positive integer $m>m_0$ and take any indices $n$ and $k$
with $n_0+1 \le n,k \le \pm$.
By (\ref{supp1}) and (\ref{supp2}), 
$$
\supp(z_k^*) \cap \supp(Tz_n) \subset \O(r(m)).
$$
Then
$$
\ymks(\ymn) = \langle P_m z_k^* , P_m T z_n \rangle 
            = z_k^*(Tz_n).
$$
Together with (\ref{tzn-zn}) and (\ref{bounded}), this gives
$$
|\ymks(\ymn) - \d_{k,n}|   
   =   |z_k^*(Tz_n-z_n)|                                     
  \le  \|z_k^*\| \cdot \frac{1}{4M}   
  \le  1/4.
$$
The proof is complete.    \endproof

\#

Observe that $\supp(\ymn)  \subset \O(r(m))$ 
and          $\supp(\ymns) \subset \O(r(m))$
for every $m > m_0$ and $n \le \pm$. 
Therefore, for every $m>m_0$ we may consider the vectors 
$\ymn$ and $\ymns$, $n \le \pm$, as elements of the space $\loo$.
 
Summarize what we have shown. 
There are positive integers $n_0$ and $m_0$ such that, 
for every integer $m>m_0$, there is a system of cardinality $(\pm-n_0)$
in $\loo$ which is:

(i)  $(1/4)$-roughly biorthogonal,

(ii) $(2M)$-bounded.

\noindent Now we establish that in this case $|\O(r(m))|$ can not
be too large.

\begin{lemma}                                                   \label{rough}
  Let $\e \in (0,1/2)$. Let $X$ be a $k$-dimensional Banach space.
  Suppose some system $(y_n,\yns)_{n \le p}$ in X
  is $\e$-roughly biorthogonal and $M$-bounded. Then
  $$
  k \ge c_1 \log p
  $$
  for some constant $c_1=c_1(\e,M)>0$.
\end{lemma}

\proof    We may assume $\|y_n\|=1$ and $\|\yns\| \le M$ for all $n$.
Then, for any non-equal indices $k$ and $n$, one has
$$
\|y_k-y_n\|   
  \ge  (y_k^*(y_k)-y_k^*(y_n)) / \|y_k^*\|     
  \ge  (1-2\e)/M  
   =:  \d.
$$
This shows that the open balls  $y_n + (\d/2) B(X)$, $n \le p$,
are pairwise disjoint and are contained in the ball 
$(1+\d/2) B(X)$. By comparing the volumes we get
$p(\d/2)^k  \le  (1+\d/2)^k$.  Hence 
$k \ge (\log(1+2/\d))^{-1} \log p$.
This completes the proof.            \endproof

Applying Lemma \ref{rough} in our situation, we obtain that 
there is a constant $c_1=c_1(M)$ such that, for every $m>m_0$,
$$
|\O(r(m))| \ge c_1 \log(\pm-n_0).
$$
Combining with (\ref{qmpm}), we have for $m>m_0$:
$$
|\O(c\l_\pm)|  \ge  c_1 \log(\pm-n_0)
$$
(clearly, we can assume that $(\l_m)_\mg$ is an increasing sequence 
of positive integers).
Since the sequence $(\pm)_\mg$ is increasing, we get
\begin{eqnarray}
  \limsup_n \frac{|\O(c\l_n)|} {\log n} > 0.            \label{limsup}
\end{eqnarray}

\noindent Now we show that there is a permutation $\pi$ on $\N$ 
such that (\ref{limsup}) fails for any constant $c$.

\begin{lemma}
  Let $f: \N \ra \R_+$ be a non-decreasing function 
  with $\lim_n f(n) = \infty$.
  Then there is a permutation $\pi$ on $\N$ with
  $$
  \lim_n \frac {|\O(cn)|} {f(n)} = 0
  $$
  for every (positive integer) constant $c$.
\end{lemma}

\proof      One can easily construct a non-decreasing 
"onto" function $\f: \N \ra \N$ such that

(i) $\lim_n \f(n) = \infty$;

(ii) $\lim_n \f(n) / f(n) = 0$;

(iii) $\f(n) \le n$ and $\f(2n) \le 2\f(n)$ for every $\ng$.

\noindent   Then we define a function $\F: \N \ra \N$:
$$
\F(n) = | \{m : \f(m) \le n \} |.
$$
Note that, for every $\mg$,
\begin{eqnarray}
  | \{n : \F(n) \le m \} |      =     \f(m).                       \label{fi}
\end{eqnarray}
Let $\G$ be a subset of $\N$ such that for every $\mg$
\begin{eqnarray}
  | \G \cap \{1,\ldots,m\} |   \le    \f(m).                    \label{gamma}
\end{eqnarray}
Finally, define a permutation $\pi$ as follows:
$$
\pi(n) = \left\{  \begin{array}{ll}
                    \F(n),     & n \not\in \G \\
                    \free(n),  & n \in \G,
                  \end{array}
         \right. 
$$
where $\free(n)$ denotes the minimal positive integer 
$k \not\in \{\pi(1), \pi(2), \ldots, \pi(n-1)\}$.
The permutation $\pi$ is well defined: 
indeed, the functions $\F(n)$ and $\free(n)$ are strictly increasing
and $\F(n) \ge n$, $\free(n) \le n$ for every $n$.
We see that
$$
\O(m) \subset \{n: \F(n) \le m\}  \cup  
  \{ \free(n): n \in \G \cap \{1,\ldots,m\} \} .
$$
By (\ref{fi}) and (\ref{gamma}),
\begin{eqnarray}
  |\O(m)| \le 2\f(m).                                           \label{omega}
\end{eqnarray}
Let a positive integer $k$ be so that $2^k \ge c$. 
By (\ref{omega}) and (iii), we have for every $\ng$:
$$
|\O(cn)|  \le  |\O(2^k n)|  \le  2\f(2^k n)  \le  2 \cdot 2^k\f(n).
$$
Thus
$$
\lim_n \frac {|\O(cn)|} {f(n)}   \le  
2^{k+1} \lim_n \frac {\f(n)} {f(n)} = 0.
$$
The proof is complete.                \endproof

\#
{\em e-mail addresses:}
paoter@mate.polimi.it, \ boris.v.novikov@univer.kharkov.ua

\end{document}